\documentclass{amsart}

\usepackage{amsfonts,txfonts,psfrag,color}
\usepackage[dvips]{graphicx}
\usepackage{epstopdf}

\begin{document}
\title{A Bound for orderings of reidemeister moves }
\author{Julian Gold}

\thanks{Research partially supported by NSF VIGRE grant DMS0636297.}

\maketitle

\begin{abstract}

We provide an upper bound on the number of ordered Reidemeister moves required to pass between two diagrams of the same link. This bound is in terms of the number of unordered Reidemeister moves required.

\end{abstract}

\vspace{10 mm}

\noindent In 1927 Kurt Reidemeister proved that any two link diagrams representing the same link may be joined by a finite sequence of Reidemeister moves. One cannot overstate the importance of this theorem to knot theory. Mathematicians like Alexander Coward \cite{ordering, upperrm}, Marc Lackenby \cite{upperrm}, Bruce Trace \cite{trace}, Joel Hass and Jeffery Lagarias \cite{HL} have all explored properties of sequences of Reidemeister moves.  

\begin{figure}[h]
\centering
\psfrag{a}[B][bl]{$\longrightarrow$}
\psfrag{b}[B][bl]{$\longleftarrow$}
\psfrag{c}[B][bl]{$\longleftrightarrow$}
\psfrag{d}[B][bl]{$\Omega_1^\uparrow$}
\psfrag{e}[B][bl]{$\Omega_1^\downarrow$}
\psfrag{f}[B][bl]{$\Omega_2^\uparrow$}
\psfrag{g}[B][bl]{$\Omega_2^\downarrow$}
\psfrag{h}[B][bl]{$\Omega_3$}
\includegraphics[scale=1.0]{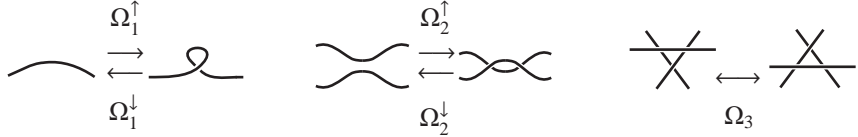}
\caption{Reidemeister Moves.} \label{F:1}
\end{figure}

\noindent In 2006, Alexander Coward showed \cite{ordering} that given any sequence of Reidemeister moves between link diagrams $D_1$ and $D_2$, it is possible to construct a new sequence ordered in the following way: first $\Omega_1^\uparrow$ moves, then $\Omega_2^\uparrow$ moves, then $\Omega_3$ moves, finally $\Omega_2^\downarrow$ moves. We present, via the following theorem, an upper bound on the number of moves required for an ordered sequence in terms of the number of moves present in any sequence of Reidemeister moves.\\

\noindent \textbf{Theorem 1. }\emph{ Let $D_1$ and $D_2$ be diagrams for the same link that are joined by a sequence of $M$ Reidemeister moves. Let $N = 6^{M+1}M$. Then there exists a sequence of no more than $\textrm{exp}^{(N)}(N)$ moves from $D_1$ to $D_2$ ordered in the following way: first $\Omega_1^\uparrow$, then $\Omega_2^\uparrow$, then $\Omega_3$, then $\Omega_2^\downarrow$ and finally $\Omega_1^\downarrow$}.\\

\noindent Here the function $\textrm{exp}$ is defined as $\textrm{exp}(x) = 2^x$ and $\textrm{exp}^{(r)}(x)$ is the function $\textrm{exp}$ iterated $r$ times on input $x$.\\

\noindent I am extraordinarily grateful to Alexander Coward for many insightful discussions and for his guidance in writing this paper. \\

\noindent We define a link diagram to be a 4-valent graph embedded in $\mathbb{R}^2$ with crossing information recorded at each vertex. All diagrams will be oriented, so that they represent oriented links. We regard two diagrams as the same if there is an ambient isotopy of $\mathbb{R}^2$ taking one diagram to the other, preserving crossing information and the orientation of each link component. To prove Theorem 1, we will adapt the methods Alexander Coward uses in \cite{ordering} and borrow the following terminology. \\

\noindent \textbf{Definition: } Let $D$ be a link diagram and suppose $c : [0,1] \to \mathbb{R}^2$ is an embedded path whose image $C$ intersects $D$ transversely at finitely many points, where $c(0) \in D$ and $c(1) \notin D$. We stipulate that no point of intersection of $D$ and $C$ is a vertex of $D$. At each such point, apart from $c(0)$, we designate whether $C$ passes over or under $D$.  \\

\noindent Let $C \times [-\epsilon,\epsilon]$ be a small neighborhood of $C$ such that $(C \times [-\epsilon,\epsilon]) \cap D = (C \cap D) \times [-\epsilon, \epsilon]$. Then define the diagram $D'$ as the 4-valent graph
\[ D \cup \partial(C \times [-\epsilon,\epsilon]) \setminus (c(0) \times(-\epsilon,\epsilon)) \]
\noindent with crossing information induced by the path $c$. We write $D \rightsquigarrow D'$ and say that $D'$ is obtained from $D$ by \emph{adding a tail along C}. Additionally, we will call $C$ the \emph{core of this tail}. We require that adding a tail to a diagram $D$ produces a diagram $D'$ where $c(D') > c(D)$. Figure 2 illustrates the construction of a tail.

\begin{figure}[h]
\centering
\psfrag{a}[B][bl]{$$}
\psfrag{b}[B][bl]{$$}
\psfrag{c}[B][bl]{$C$}
\psfrag{d}[B][bl]{$\Omega_2^\downarrow$}
\psfrag{e}[B][bl]{$\Omega_3$}
\includegraphics[scale=1.5]{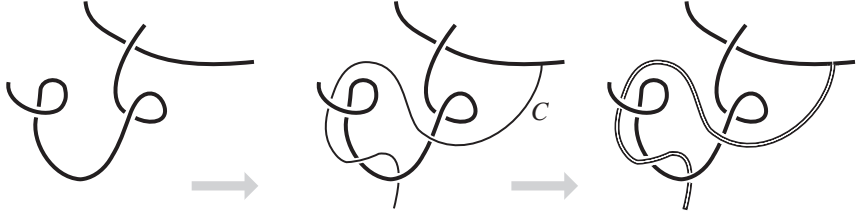}
\caption{Adding a tail.} \label{F:2}
\end{figure}

\noindent \textbf{Definition: } Suppose $D_1 \rightsquigarrow D_2$ via some path $c : [0,1] \to \mathbb{R}^2$. Suppose additionally that $c(1)$ lies in a small neighborhood of some crossing $\chi$ of $D_1$. Let $D_3$ be as in Figure 3, a diagram obtained from $D_2$ by performing two $\Omega_2^\uparrow$ moves followed by one $\Omega_3$ move:\\

\begin{figure}[h]
\centering
\psfrag{a}[B][bl]{\large{$D_1$}}
\psfrag{b}[B][bl]{\large{$D_2$}}
\psfrag{c}[B][bl]{\large{$D_3$}}
\psfrag{d}[B][bl]{$\xrightarrow{\Omega_2^\uparrow , \Omega_2^\uparrow, \Omega_3}$}
\psfrag{o}[B][bl]{\Large{$\rightsquigarrow$}}
\psfrag{e}[B][bl]{$\chi$}
\includegraphics[scale=1.5]{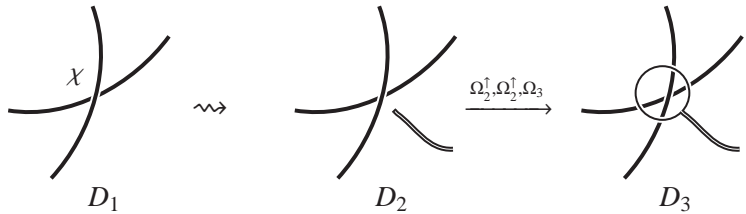}
\caption{Adding a lollipop.} \label{F:3}
\end{figure}

\noindent We say $D_3$ is obtained from $D_1$ by \emph{adding a lollipop} and write $D_1  \circleright D_3$. The \emph{lollipop} itself is defined as $\overline{D_3 \setminus D_1}$. The \emph{tail part of the lollipop} is $\overline{(D_3 \cap D_2) \setminus D_1}$, and the closure of the rest of the lollipop is the \emph{circle part of the lollipop}. We say that the lollipop is \emph{centered at $\chi$}.\\

\noindent We think of a sequence $\mathcal{S}$ of Reidemeister moves, tails and lollipops between link diagrams $L_1$ and $L_2$ in the following way:
\[ \mathcal{S} : L_1 = D_0 \xrightarrow{a_1} D_1 \xrightarrow{a_2} \dots \xrightarrow{a_n} D_n = L_2 \]
\noindent Here each $a_i$ is a Reidemeister move, a tail or a lollipop. A tail or lollipop may be added from $D_i$ to $D_{i+1}$ (eg $D_i \rightsquigarrow D_{i+1}$) or from $D_{i+1}$ to $D_i$ (eg $D_i \leftsquigarrow D_{i+1}$). We say the \emph{length} of $\mathcal{S}$ is $n$.  The intermediate link diagrams $D_i$ are often omitted from the figures in this paper for clarity, but are implicit in any sequence.\\
  
\noindent If a link diagram $D_2$ is reached from $D_1$ by a sequence of $\Omega_2^\uparrow$ moves of length $n$, we write $D_1 \twoheadrightarrow^{n} D_2$.\\

\noindent The following lemma allows us to take a sequence $\mathcal{S}$ and produce a sequence $\mathcal{S}'$ with one less $\Omega_3$ move.\\

\noindent $\textbf{Lemma 2. }$ \emph{Let $D_1$ and $D_2$ be link diagrams such that $D_1 \xrightarrow{\Omega_3} D_2$. Then there exists a diagram $D_3$ such that $D_1 \twoheadrightarrow^{2} D_3$ and $D_2 \circleright D_3$. }

\begin{proof}
\hspace{20mm}
\begin{figure}[h]
\centering
\psfrag{a}[B][bl]{$D_1$}
\psfrag{b}[B][bl]{$D_2$}
\psfrag{c}[B][bl]{$D_3$}
\psfrag{d}[B][bl]{\footnotesize{$2$}}
\psfrag{e}[B][bl]{\footnotesize{$\Omega_3$}}
\includegraphics[scale=1.2]{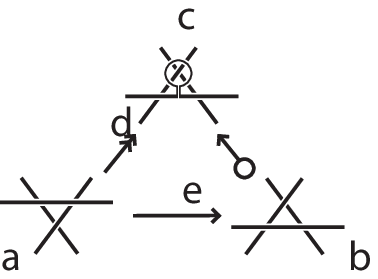}
% \caption{} \label{F:4}
\end{figure}

\end{proof}

\noindent If an $\Omega_3$ move occurs in a sequence of Reidemeister moves, tails and lollipops $\mathcal{S}$, we can apply Lemma 2 to $\mathcal{S}$ to get a new sequence $\mathcal{S}'$:
\[ \mathcal{S} : A \to \dots \to B \xrightarrow{\Omega_3} C \to \dots \to D \]
\[ \mathcal{S}' :A \to \dots \to B \xrightarrow{\Omega_2^\uparrow} B'  \xrightarrow{\Omega_2^\uparrow} B'' \circleleft C \to \dots \to D \]
When we apply Lemma 2 to construct $\mathcal{S}'$ from $\mathcal{S}$, we call this \emph{capping the $\Omega_3$ move} from $B$ to $C$. \\

% Added the following sentence, introducing Prop. 3:
\noindent The following proposition and its corollary will also allow us to build new sequences from old ones in a useful way.\\

\noindent $\textbf{Proposition 3. }$ \emph{ Suppose $D_1 \rightsquigarrow D_1'$ (or $D_1 \circleright D_1'$) and also that $D_1 \twoheadrightarrow^{1} D_2$. Then there exists a diagram $D_2'$ such that $D_2 \rightsquigarrow D_2'$ ($D_2 \circleright D_2'$ respectively) and $D_1'  \twoheadrightarrow^\alpha D_2'$, where }
\begin{equation*}  c(D_2') - c(D_2) \leq 2( c(D_1') - c(D_1) ) \tag{\textbf{A}}\end{equation*}
\noindent \emph{and}
\begin{equation*}   \alpha \leq c(D_1') - c(D_1). \tag{\textbf{B}}\end{equation*}

\begin{figure}[h]
\centering

\psfrag{d}[B][bl]{$D_1'$}
\psfrag{e}[B][bl]{$D_2'$}
\psfrag{i}[B][bl]{$D_1$}
\psfrag{j}[B][bl]{$D_2$}
\psfrag{k}[B][bl]{\footnotesize{$\alpha$}}
\psfrag{o}[B][bl]{\footnotesize{$1$}}

\includegraphics[scale=.8]{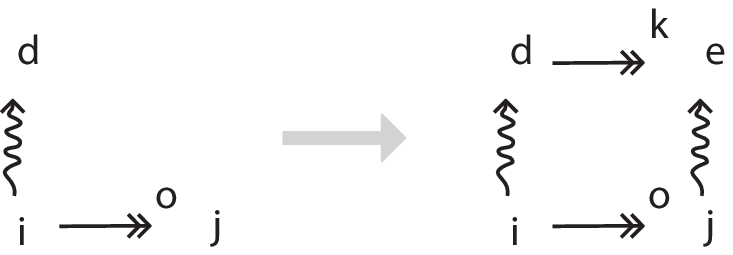}
% \caption{An illustration of Proposition 3, in the case $D_1 \rightsquigarrow D_1'$.} \label{F:8}
\end{figure}

\begin{proof} 
The diagram $D_2$ is obtained from $D_1$ by a single $\Omega_2^\uparrow$ move which takes place over two (possibly non-distinct) edges $e_1$ and $e_2$ of $D_1$. Pick points $p_1$ and $p_2$ on $e_1$ and $e_2$ respectively, so that $p_1$ and $p_2$ are disjoint from a neighborhood of the tail $D_1 \rightsquigarrow D_1'$. We can perform the $\Omega_2^\uparrow$ move from $D_1$ to $D_2$ by adding a tail along a path $P$, which starts at $p_1$ and ends slightly beyond $p_2$.

\begin{figure}[h]
\centering
\psfrag{a}[B][bl]{$e_1$}
\psfrag{b}[B][bl]{$e_2$}
\psfrag{c}[B][bl]{$p_1$}
\psfrag{d}[B][bl]{$p_2$}
\psfrag{e}[B][bl]{$P$}
\includegraphics[scale=1.5]{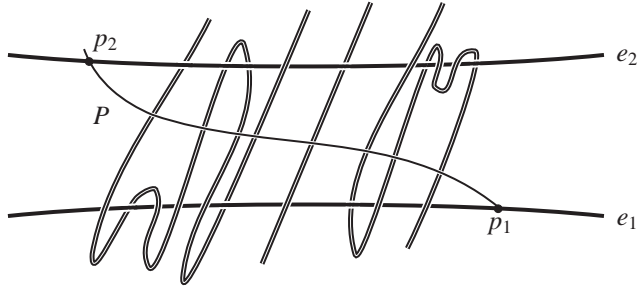}
\caption{Constructing $D_2'$ by adding a tail along $P$.} \label{F:5}
\end{figure}

\noindent Diagram $D_1'$ contains the points $p_1$ and $p_2$. We may arrange that the intersection of $P$ with the tail $D_1 \rightsquigarrow D_1'$ contains at most $2 \lfloor \frac{c(D_1') - c(D_1)}{4} \rfloor $ points. Figure 4 depicts such an arrangement. Adding a tail along $P$, we construct a diagram $D_2'$ with 
\[ c(D_2') - c(D_1') \leq 4 \lfloor \frac{c(D_1') - c(D_1)}{4} \rfloor +2. \]
\noindent Hence
\[ c(D_2') - c(D_1') \leq c(D_1') - c(D_1) +2. \]
\noindent We note that $c(D_1') - c(D_1) +2 \leq 2(c(D_1') - c(D_1))$, because adding a tail to a diagram must raise its crossing number by at least two. This implies the desired bound on $\alpha$.
Also \[ c(D_2') - c(D_1') \leq c(D_1') - c(D_1) +2 \]
implies, by adding $c(D_1')$ to both sides and subtracting $c(D_2)$, that
 \[ c(D_2') - c(D_2) \leq 2c(D_1') - c(D_1) +2 - c(D_2). \] 
Using $c(D_2) = c(D_1) + 2$ we get
\[ c(D_2') - c(D_2) \leq 2c(D_1') - 2c(D_1). \] 
% Added: "...by adding $c(D_1')$ to and subtracting $c(D_2)$ from both sides, that..."

\noindent In the case that $D_1 \circleright D_1'$, choose $p_1$ and $p_2$ to be outside the circle part of the lollipop, and the above considerations go through.

\end{proof}

% Added the following sentence, introducing Corollary 4:
\noindent Corollary 4 is a natural generalization of Proposition 3.\\

%  Changed alpha to beta in the statement of Corollary 4
\noindent $\textbf{Corollary 4. }$  \emph{Suppose $D_1 \rightsquigarrow D_1'$ (or $D_1 \circleright D_1')$ and also that $D_1 \twoheadrightarrow^{n} D_2$. Then there exists a diagram $D_2'$ such that $D_2 \rightsquigarrow D_2'$ ($D_1 \circleright D_1' respectively$) and $D_1'  \twoheadrightarrow^\beta D_2'$, where }
\[ \beta \leq 2^n (c(D_1') - c(D_1) ). \]

\begin{figure}[h]
\centering

\psfrag{d}[B][bl]{$D_1'$}
\psfrag{e}[B][bl]{$D_2'$}
\psfrag{i}[B][bl]{$D_1$}
\psfrag{j}[B][bl]{$D_2$}
\psfrag{k}[B][bl]{\footnotesize{$\beta$}}
\psfrag{o}[B][bl]{\footnotesize{$n$}}

\includegraphics[scale=.8]{comm1.eps}
% \caption{Corollary 4, in the case $D_1 \rightsquigarrow D_1'$.} \label{F:8}
\end{figure}

\begin{proof}
Let $D_1, D_2$ and $D_1'$ be as in the statement of the theorem. We work in the case $D_1 \rightsquigarrow D_1'$, but the proof for lollipops is identical. Let $\mathcal{S}$ be the sequence of $\Omega_2^\uparrow$ moves of length $n$ from $D_1$ to $D_2$,
\[ \mathcal{S}: D_1 = E_0 \twoheadrightarrow^1 E_1 \twoheadrightarrow^1 \dots \twoheadrightarrow^1 E_n = D_2, \]

% Changed $\beta$ to $\beta_0$ in the following paragraph:
% Added "..to build a diagram $E_2'$.."
% Changed some sentence structure of paragraph
\noindent and let $E_0' = D_1'$. We use Proposition 3 to construct a diagram $E_1'$ such that $E_1 \rightsquigarrow E_1'$ and $E_0' \twoheadrightarrow^{\beta_0} E_1'$, where $\beta_0 \leq c(E_0') - c(E_0)$. Apply Proposition 3 again to the triple $(E_1,E_1',E_2)$ to build a diagram $E_2'$. Repeat this application to construct the diagrams $E_2'$ through $E_n'$, as below.

\begin{figure}[h]
\centering
\psfrag{a}[B][bl]{$D_1' = E_0'$}
\psfrag{b}[B][bl]{$E_1'$}
\psfrag{c}[B][bl]{$\dots$}
\psfrag{d}[B][bl]{$E_{n-2}'$}
\psfrag{e}[B][bl]{$E_{n-1}'$}
\psfrag{f}[B][bl]{$E_n' = D_2'$}
\psfrag{g}[B][bl]{$D_1 = E_0$}
\psfrag{h}[B][bl]{$E_1$}
\psfrag{i}[B][bl]{$E_{n-2}$}
\psfrag{j}[B][bl]{$E_{n-1}$}
\psfrag{k}[B][bl]{$E_n = D_2$}
\psfrag{o}[B][bl]{\footnotesize{$1$}}
\psfrag{t}[B][bl]{\footnotesize{$\beta_0$}}
\psfrag{u}[B][bl]{\footnotesize{$\beta_1$}}
\psfrag{r}[B][bl]{\footnotesize{$\beta_{n-2}$}}
\psfrag{s}[B][bl]{\footnotesize{$\beta_{n-1}$}}
\psfrag{v}[B][bl]{$\beta$}

\includegraphics[scale=.8]{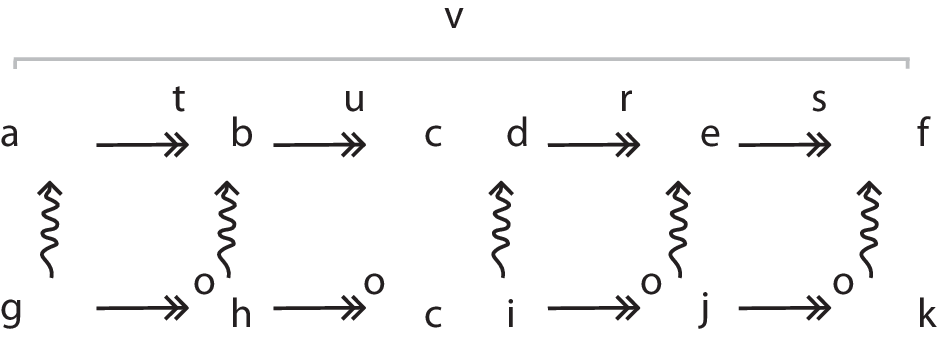}
% \caption{} \label{F:8}
\end{figure}

\noindent Proposition 3 ($\textbf{B}$) gives us that $\beta_i \leq c(E_i') - c(E_i)$, while proposition 3 ($\textbf{A}$) tells us $c(E_i') - c(E_i) \leq 2^i(c(E_0') - c(E_0))$. The sequence of $\Omega_2$ moves from $E_0'$ to $E_n'$ has length $\beta$, where $\beta = \sum_{i=0}^{n-1} \beta_i$. Hence, 
\[ \beta \leq (2^n - 1) (c(E_0') - c(E_0)). \]
\noindent Take $D_2' = E_n'$ and a larger bound on $\beta$ to complete the proof. 

\end{proof}

% -----------------------------------------------

\noindent Theorem 5 uses Lemma 2, Proposition 3 and Corollary 4 to begin building an ordered sequence from an unordered sequence.\\

\noindent \textbf{Theorem 5. } \emph{ Let $D_2$ be a link diagram obtained from $D_1$ via a sequence of $\Omega_2$ and $\Omega_3$ moves of length $M$. Then there exists a diagram $D_3$ such that $D_1 \twoheadrightarrow^{\gamma} D_3$ and $D_3$ is obtained from $D_2$ by adding a sequence no more than $M$ tails and lollipops. Further,}
\[ \gamma \leq \text{exp}^{(M)}(6M).\]

\begin{proof}

Consider a sequence $\mathcal{A}$ of $\Omega_2$ and $\Omega_3$ moves of length $M$ from $D_1$ to $D_2$, $\alpha_3$ of which are $\Omega_3$:
\[ \mathcal{A} : D_1 = A_0 \to A_1 \to \dots \to A_M = D_2 \]
\noindent Using Lemma 2, cap every $\Omega_3$ move to build a new sequence $\mathcal{E}_1$ with no $\Omega_3$ moves:
\[ \mathcal{E}_1 : D_1 = E_0 \to E_1 \to \dots \to E_{M+2\alpha_3} = D_2 \]

\begin{figure}[h]
\centering
\psfrag{a}[B][bl]{$\dots$}
\psfrag{b}[B][bl]{$A_0$}
\psfrag{c}[B][bl]{$A_M$}
\psfrag{d}[B][bl]{$E_0$}
\psfrag{e}[B][bl]{$E_{M + 2\alpha_3}$}
\psfrag{f}[B][bl]{$A_1$}
\psfrag{g}[B][bl]{$A_2$}
\psfrag{h}[B][bl]{$E_1$}
\psfrag{i}[B][bl]{$E_3$}
\psfrag{j}[B][bl]{$E_4$}

\includegraphics[scale=1.0]{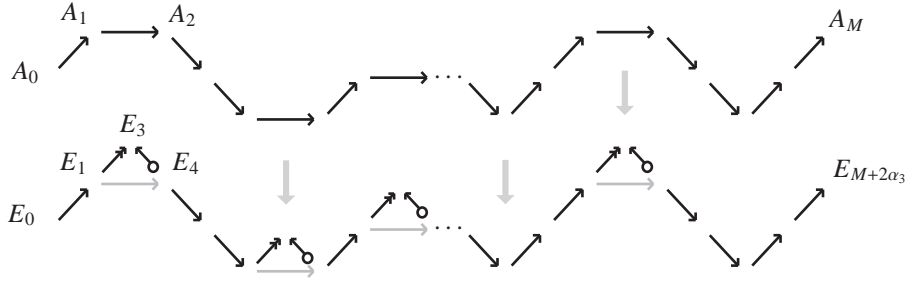}
\caption{Constructing $\mathcal{E}_1$ from $\mathcal{A}$.} \label{F:11}
\end{figure}

\noindent  If $E_i \xrightarrow{\Omega_2^\downarrow} E_{i+1}$, we instead write $E_i \leftsquigarrow E_{i+1}$, because a $\Omega_2^\uparrow$ move may be performed by adding a tail. Define a \emph{local minimum} of $\mathcal{E}_1$ to be a diagram $E_i$ such that 
\[ E_{i-1} \leftsquigarrow E_i  \xrightarrow{\Omega_2^\uparrow} E_{i+1} \text{ \hspace{10mm} or \hspace{10mm} }  E_{i-1} \circleleft E_i  \xrightarrow{\Omega_2^\uparrow} E_{i+1}.\]

\noindent Let $E_x \in \{ E_1, \dots, E_{M+2\alpha_3-1} \}$ be the last local minimum appearing in $\mathcal{E}_1$. Let $r_1$ be the number of consecutive $\Omega_2^\uparrow$ moves in $\mathcal{E}_1$ to the right of $E_x$. Let $\ell_1$ be the number of consecutive $\Omega_2^\uparrow$ moves in $\mathcal{E}_1$ to the left of $E_{x-1}$. \\

\begin{figure}[h]
\centering
\psfrag{a}[B][bl]{$\dots$}
\psfrag{b}[B][bl]{$D_1$}
\psfrag{c}[B][bl]{$D_2$}
\psfrag{d}[B][bl]{$E_0$}
\psfrag{e}[B][bl]{$E_{x + r_1}$}
\psfrag{f}[B][bl]{$E_x$}
\psfrag{g}[B][bl]{$F$}
\psfrag{h}[B][bl]{$$}
\psfrag{i}[B][bl]{$r_2$}
\psfrag{j}[B][bl]{$E_{x-1}$}
\psfrag{k}[B][bl]{$$}
\psfrag{z}[B][bl]{$r_1$}

\includegraphics[scale=1.0]{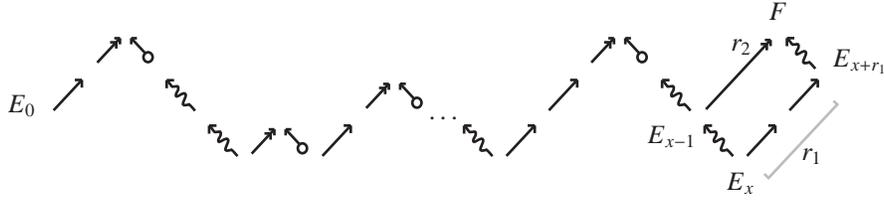}
\caption{Constructing $F$. In this case, $E_{x+r_1} = E_{M+2\alpha_3}$.} \label{F:11}
\end{figure}

\noindent Apply Corollary 4 to the triple $(E_{x-1}, E_x, E_{x+ r_1})$ to build a diagram $F$, where $E_{x-1} \twoheadrightarrow^{r_2'} F$ and where $E_{x+r_1} \circleright F$ if $E_x \circleright E_{x-1}$ or $E_{x +z_1} \rightsquigarrow F$ if $E_x \rightsquigarrow E_{x-1}$. Corollary 4 tells us $r_2' \leq 4 \cdot 2^{r_1}$, in the worst case that $E_x \circleright E_{x-1}$. Figure 6 depicts the construction of $F$.\\

\noindent Define $\mathcal{E}_2$ to be the following sequence:
\[ \mathcal{E}_2 : D_1 = E_0 \to E_1 \to \dots \to E_{x-1} \to \dots \to F \to E_{x+r_1} \to \dots \to E_{M+2\alpha_3} \]

\noindent Then $\mathcal{E}_2$ is a sequence of diagrams with $r_2$ consecutive $\Omega_2^\uparrow$ moves to the right of its last local minimum, with $r_2$ bounded by
\[ r_2 \leq 4 \cdot 2^{r_1} + \ell_1.\]
Hence 
\[r_2 \leq 2^{r_1 + 2 + \ell_1}\]
\noindent In general let $\mathcal{E}_k$ be a sequence with $r_k$ the number of $\Omega_2^\uparrow$ moves to the right of the last local minimum of $\mathcal{E}_k$. Let $\ell_k$ be the number of consecutive $\Omega_2^\uparrow$ moves preceding the diagram to the immediate left of the last local minimum of $\mathcal{E}_k$. Given the pair $(\mathcal{E}_k, r_k)$, we may apply Corollary 4 as above to produce a pair $(\mathcal{E}_{k+1}, r_{k+1})$ satisfying
\[ r_{k+1} \leq 2^{r_k + 2 + \ell_k}. \]
\noindent Inductively, 
\[ r_{k+1} \leq \text{exp}^{(k)}\left(r_1 + 2k + \sum_{i=1}^k \ell_i \right). \]

\begin{figure}[h]
\centering
\psfrag{a}[B][bl]{$\dots$}
\psfrag{b}[B][bl]{$D_1$}
\psfrag{c}[B][bl]{$D_2$}
\psfrag{d}[B][bl]{$E_0$}
\psfrag{e}[B][bl]{$E_{x + z}$}
\psfrag{f}[B][bl]{$E_x$}
\psfrag{g}[B][bl]{$F_1$}
\psfrag{h}[B][bl]{$$}
\psfrag{i}[B][bl]{$D_3$}
\psfrag{j}[B][bl]{$$}
\psfrag{k}[B][bl]{$$}
\psfrag{l}[B][bl]{$$}
\psfrag{m}[B][bl]{$$}
\psfrag{x}[B][bl]{\Large{$r_1$}}
\psfrag{z}[B][bl]{\Large{$r_{\pi}$}}

\includegraphics[scale=1.0]{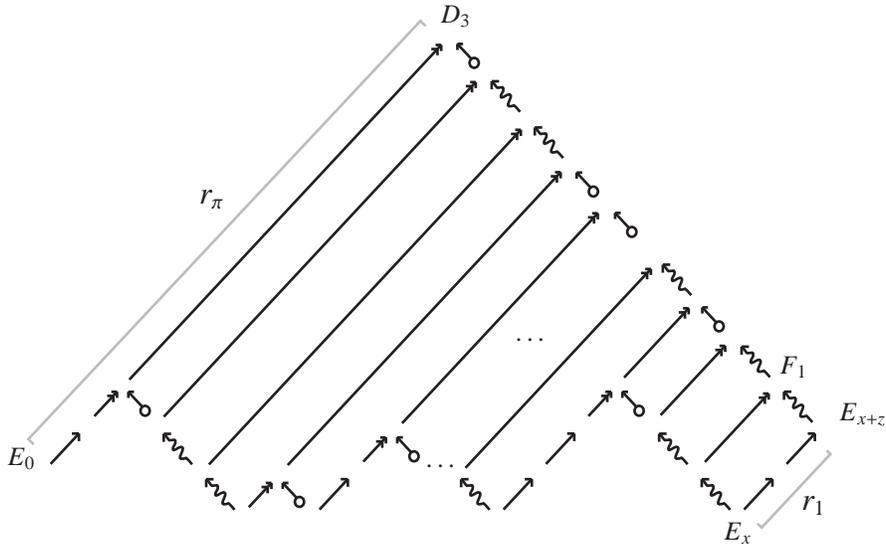}
\caption{Repeatedly applying Corollary 4 to build $D_3$.} \label{F:11}
\end{figure}

\noindent Iterate the constructions of the $(\mathcal{E}_k,r_k)$ until we produce a sequence $\mathcal{E}_\pi$ with no local minima and with $r_\pi$ consecutive $\Omega_2^\uparrow$ moves following $E_0$. The number of times we apply Corollary 4 to construct $\mathcal{E}_\pi$ from $\mathcal{E}_1$ is exactly the number of tails and lollipops in $\mathcal{E}_1$, which is less than or equal to $M$. So $\pi \leq M+1$, and via our above formula,
\[ r_\pi \leq \text{exp}^{(\pi-1)}\left(r_1 + 2(\pi-1) + \sum_{i=1}^{\pi-1} \ell_i \right). \]
\noindent We note that $r_1 \leq M$ and $\sum_{i=1}^{\pi-1} \ell_i \leq M +2\alpha_3 \leq 3M$. Substituting, we get
\[ r_\pi \leq \text{exp}^{(M)}(6M). \]
\noindent There are $r_\pi$ moves of type $\Omega_2^\uparrow$ following $D_1 = E_0$ in $\mathcal{E}_\pi$, so let $D_3$ be the diagram obtained by performing these moves on $D_1$. Because $D_3$ is obtained from $E_{M+2\alpha_3} = D_2$ by at most $M$ tails and lollipops, Theorem 5 holds.

\end{proof}

% -----------------------------------------------

\noindent The following theorem allows us to construct an ordered sequence of $\Omega_2$ and $\Omega_3$ moves from the tails and lollipops arising in Theorem 5.\\

\noindent $\textbf{Theorem 6. }$ \emph{Suppose $D_2$ is obtained from $D_1$ by adding a sequence $\mathcal{T}$ of tails and lollipops of length $S$:}
\[ \mathcal{T} : D_1 = T_0 \xrightarrow{a_1} T_1 \xrightarrow{a_2} \dots \xrightarrow{a_S} T_S = D_2 \]
\noindent \emph{where either $T_i \rightsquigarrow T_{i+1}$ or $T_i \circleright T_{i+1}$. Then there exists a diagram $D_3$ obtained from $D_2$ by a sequence of $\Omega_2^\uparrow$ moves of length no more than $\frac{S}{2}(c(D_2) - c(D_1)) + 2S$, followed by a sequence of $\Omega_3$ moves of length no more than $S$. Additionally $D_1$ is obtained from $D_3$ by a sequence of $\Omega_2^\downarrow$ moves of length at most $\frac{S+1}{2}(c(D_2) - c(D_1)) + 2S$.}

\begin{proof}

\noindent Consider a crossing $\chi$ of the diagram $D_2$ about which the circle part of a lollipop in $\mathcal{T}$ is centered. There may be multiple lollipops (suppose there are $k$) centered at $\chi$, so consider a point $p_k$ on the outermost one. Let $q$ be a point in a small enough neighborhood of $\chi$ such that a straight line segment from $q$ to $\chi$ does not intersect $D_2$ except at $\chi$. \\

\noindent Consider a path $c : [0,1] \to \mathbb{R}^2$ such that $c(0) = p_k$ and $c(1) = q$. Choose $c$ in such a way that its image $C$ intersects each concentric lollipop at only one point. The point of intersection of $C$ and the $i$th concentric lollipop is denoted $p_i$. Let $\delta_k = 0$ and let $\delta_{k-1} < \delta_{k-2} < \dots < \delta_1$ be real numbers in $(0,1)$ such that $c(\delta_i) = p_i$. \\

\noindent Via the argument used in the proof of Proposition 3, we also choose $c$ so that $C \cap D_2$ consists of no more than $2 \lfloor \frac{c(D_2) - c(D_1)}{4} \rfloor$ points, excluding the points $p_1$ through $p_k$. 

\begin{figure}[h]
\centering
\psfrag{a}[B][bl]{$\longrightarrow$}
\psfrag{b}[B][bl]{$\chi$}
\psfrag{c}[B][bl]{$p_k$}
\psfrag{d}[B][bl]{$q$}
\psfrag{e}[B][bl]{$q_2$}
\psfrag{f}[B][bl]{$p_1$}
\psfrag{g}[B][bl]{\large{$D_2 \cup C$}}
\psfrag{j}[B][bl]{\large{$E$}}

\includegraphics[scale=1.5]{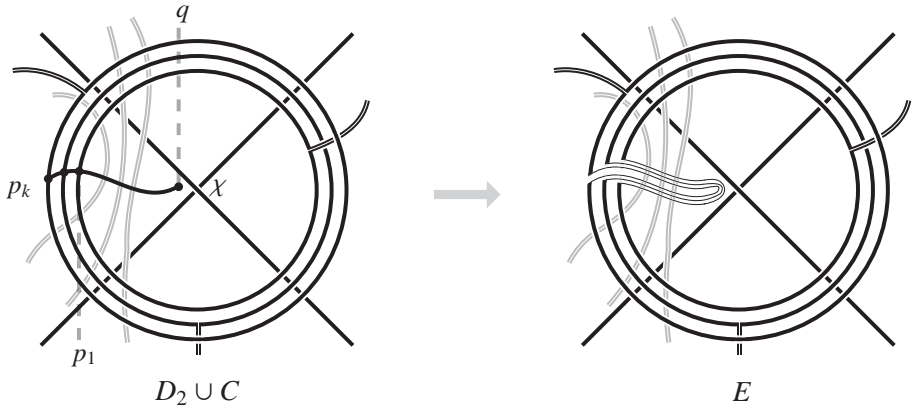}
\caption{Adding concentric tails at the crossing $\chi$.} \label{concentric}
\end{figure}

\noindent Add a tail along the path $c|_{[\delta_1,1]}$ to construct a diagram $E_1$ from $D_2$, where $c(E_1) - c(D_2) \leq c(D_2) - c(D_1)$. Perturb this tail slightly, so that it is closer to the crossing $\chi$, and now add a second tail disjoint from the first tail along the path $c|_{[\delta_2,1]}$. This second tail introduces no more than $c(D_2) - c(D_1)$ crossings. \\

\noindent Repeating this process of perturbing and adding tails along $c|_{[\delta_i,1]}$ for all $i \in [1, \dots, k]$, we produce a diagram $E_k$ where $c(E_k) - c(D_2) \leq k(c(D_2) - c(D_1))$. To build a diagram $E$, add nested tails in the same way for every crossing of $D_2$ that is the center of some lollipop, so that $c(E) - c(D_2) \leq S(c(D_2)-c(D_1))$. Then $E$ may be obtained from $D_2$ by a sequence of $\Omega_2^\uparrow$ moves of length at most $\frac{S}{2}(c(D_2)-c(D_1))$. \\

\noindent Now construct the diagram $E'$ from $E$ by doing the following for each crossing: If there are $k$ concentric circles centered at a crossing $\chi$, perform $2k$ type $\Omega_2^\uparrow$ moves, forking the previously constructed tails over the edges of the crossing $\chi$, as Figure 9 illustrates.

\begin{figure}[h]
\centering
\psfrag{a}[B][bl]{$\longrightarrow$}
\psfrag{b}[B][bl]{\large{$E$}}
\psfrag{c}[B][bl]{\large{$E'$}}
\includegraphics[scale=1.5]{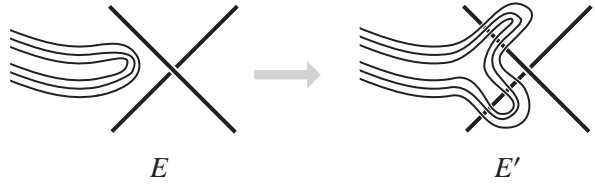}
\caption{Perform $2k$ type $\Omega_2^\uparrow$ moves, so that each tail `forks' over the crossing.} \label{forker}
\end{figure}

\noindent The diagram $E'$ may be reached from $D_2$ via a sequence of $\Omega_2^\uparrow$ moves with length at most $\frac{S}{2}(c(D_2) - c(D_1)) + 2S$. Finally, construct the diagram $D_3$ by performing at most $S$ moves of type $\Omega_3$, as in Figure 10.

\begin{figure}[h]
\centering
\psfrag{a}[B][bl]{$\longrightarrow$}
\psfrag{d}[B][bl]{\large{$E'$}}
\psfrag{e}[B][bl]{\large{$D_3$}}
\includegraphics[scale=1.5]{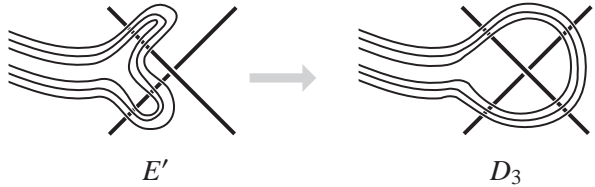}
\caption{Performing $\Omega_3$ moves to pass from $E'$ to $D_3$.} \label{F:13}
\end{figure}

\noindent We may now pass from $D_3$ to $D_1$ by performing $\Omega_2^\downarrow$ moves as follows. Each tail and lollipop of $\mathcal{T}$ in $D_2$ is still present in $D_3$, with the circle parts of each lollipop modified. We remove them one at a time starting with the last tail or lollipop $a_S$ in the sequence. If $a_S$ is a lollipop, it now has the form depicted by Figure 11 in $D_3$, and may be removed by $\Omega_2^\downarrow$ moves. If $a_S$ is a tail, it may likewise be removed by $\Omega_2^\downarrow$ moves. We continue to remove tails and lollipops in the reverse order they are added in $\mathcal{T}$ until we obtain $D_1$.\\

% First remove the last tail or lollipop $a_S$ added in $\mathcal{T}$. If $a_S$ is a lollipop, it has the form depicted in Figure 11. \\

% If the last tail or lollipop $a_S$ present in $\mathcal{T}$ is a tail, it may be removed by $\Omega_2^\downarrow$ moves. If $a_S$ is a lollipop, it has the form depicted by Figure 11 in $D_3$, and may also be removed by a sequence of $\Omega_2^\downarrow$ moves. \\

% \noindent Observe that we may pass from $D_3$ to $D_1$ by performing $\Omega_2^\downarrow$ moves. Any tail $T_i \rightsquigarrow T_{i+1}$ in $\mathcal{T}$ remains a tail in $D_3$ and can be removed by a sequence of $\Omega_2^\downarrow$ moves. In $D_3$, the circle part of a lollipop $T_j \circleright T_{j+1}$ in $\mathcal{T}$ is as in Figure 11 and may also be removed by a sequence of $\Omega_2^\downarrow$ moves. \\
\begin{figure}[h]
\centering
\includegraphics[scale=1.5]{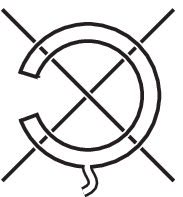}
\caption{} \label{concentric}
\end{figure}

\noindent Because $c(D_3) = c(E')$, we know $c(D_3) - c(D_1)$ is exactly $c(E') - c(D_2) + c(D_2) - c(D_1)$, which is at most $S(c(D_2) - c(D_1)) + 4S + c(D_2) - c(D_1)$. Simplifying and then halving this bound gives us the number of $\Omega_2^\downarrow$ moves from $D_3$ to $D_1$. 

\end{proof}

\noindent We consolidate previous results into Theorem 7, which is a special case of Theorem 1.\\

\noindent $\textbf{Theorem 7. }$ \emph{Let $D_2$ be a link diagram obtained from $D_1$ by a sequence of $\Omega_2$ and $\Omega_3$ moves of length $M$. Then there is a sequence of at most} $\text{exp}^{(2M)}(6M)$ \emph{Reidemeister moves from $D_1$ to $D_2$ ordered in the following way: first $\Omega_2^\uparrow$ moves, then $\Omega_3$ moves and finally $\Omega_2^\downarrow$ moves. }

\begin{proof}
Given $D_1$ and $D_2$, construct a diagram $D_3$ using Theorem 5, where $D_3$ is obtained from $D_1$ by no more than $\text{exp}^{(M)}(6M)$ type $\Omega_2^\uparrow$ moves. Additionally, $D_3$ is obtained from $D_2$ by no more than $M$ tails and lollipops. Note that $c(D_3) - c(D_2) \leq 2\cdot\text{exp}^{(M)}(6M) + 2M$.\\

% Changed $2M$ to $M$
\noindent From $D_2$ and $D_3$, apply Theorem 6 to construct a diagram $D_4$ with the following properties. There is a sequence of $\Omega_2^\uparrow$ moves whose length is no more than $M\cdot\text{exp}^{(M)}(6M) + M^2 + 2M$, followed by a sequence of $\Omega_3$ moves of length no more than $M$ from $D_3$ to $D_4$. There is also a sequence of $\Omega_2^\downarrow$ moves whose length is at most $(M+1)\cdot\text{exp}^{(M)}(6M) + M^2 + 3M$ from  $D_4$ to $D_2$. \\

% Changed 2M + 4 to 2M + 2
\noindent Following the sequences of moves constructed from $D_1$ to $D_3$, then to $D_4$ and finally to $D_2$, we have a sequence of no more than $(2M+2)\cdot\text{exp}^{(M)}(6M) + (2M+6)M$ Reidemeister moves ordered as desired. For $M \geq 1$, $\text{exp}^{(2M)}(6M) \geq (2M+2)\cdot\text{exp}^{(M)}(6M) + (2M+6)M$.

\end{proof}

% -----------------------------------------------

\noindent Before considering the more general case of an arbitrary sequence of $M$ Reidemeister moves, we need two lemmas relating to $\Omega_1$ moves. These lemmas allow us to take a sequence of Reidemeister moves and build a new sequence in which the $\Omega_1$ moves occur only at the beginning and end.\\

\noindent \textbf{Lemma 8. } \emph{Let A, B and C be link diagrams such that}
\[ A \xrightarrow{\Omega} B \xrightarrow{\Omega_1^\uparrow} C\]
\noindent \emph{where $\Omega$ is an arbitrary $\Omega_2$ or $\Omega_3$ move. Then there exists a diagram $B'$ which may be obtained from $A$ by a single $\Omega_1^\uparrow$ move, and where $C$ is obtained from $B'$ by no more than six $\Omega_2$ or $\Omega_3$ moves. Additionally, if instead $\Omega = \Omega_1^\downarrow$, there is a diagram $B'$ such that}
\[ A \xrightarrow{\Omega_1^\uparrow} B' \xrightarrow{\Omega_1^\downarrow} C. \]

\noindent \textbf{Lemma 9. } \emph{Let $A$, $B$ and $C$ be link diagrams such that $A \xrightarrow{\Omega_1^\downarrow} B \xrightarrow{\Omega} C$, where $\Omega$ is an $\Omega_2$ or $\Omega_3$ move. Then there exists a diagram $B'$ such that $B'$ is obtained form $A$ by no more than six $\Omega_2$ or $\Omega_3$ moves and where $C$ may be obtained from $B'$ by a single $\Omega_1^\downarrow$ move.}\\

\noindent The proofs of Lemma 8 and Lemma 9 are left to be verified by the reader, and Corollary 10 is a rapid consequence of these lemmas. \\

\noindent \textbf{Corollary 10. } \emph{Let $D_2$ be obtained from $D_1$ by an arbitrary sequence of $M$ Reidemeister moves, $\alpha$ of which are $\Omega_1^\uparrow$ and $\beta$ of which are $\Omega_1^\downarrow$. Then there exist diagrams $D_1'$ and $D_2'$ such that $D_1'$ is obtained from $D_1$ by $\alpha$ type $\Omega_1^\uparrow$ moves and $D_2$ is obtained from $D_2'$ by $\beta$ type $\Omega_1^\downarrow$ moves. Additionally, $D_2'$ is obtained from $D_1'$ by no more than $6^M M$ Reidemeister moves of type $\Omega_2$ and $\Omega_3$.}\\

\noindent \emph{Proof of Theorem 1.}  Begin with an arbitrary sequence of $M$ Reidemeister moves from diagram $D_1$ to diagram $D_2$, $\alpha$ of which are $\Omega_1^\uparrow$ and $\beta$ of which are $\Omega_1^\downarrow$. Construct $D_1'$ and $D_2'$ as in Corollary 10. Then apply Theorem 7 to the sequence of $\Omega_2$ and $\Omega_3$ moves from $D_1'$ to $D_2'$ to obtain a sorted sequence of Reidemeister moves from $D_1$ to $D_2$ of length at most 
\[ \text{exp}^{(6^MM)}(6\cdot 6^MM) + \alpha + \beta \leq \text{exp}^{(6^{M+1}M)}(6^{M+1}M). \]
\qed

% -----------------------------------------------

\bibliographystyle{amsplain}
\bibliography{rmrefs}
% \section*{\Large{c(L)}}
% \label{sec:c(L)}
% \input {cL/cL}

% \section*{\Large{c(P)}}
% \label{sec:c(P)}
% \input {cP/cP}
 
% \section*{\Large{Putting it together}}
% \label{sec:Putting it together}
% \input {Tog/Tog}

% \section*{\Large{$\Omega_1$ Moves}}
% \label{sec:Omega1 Moves}
% \input{o1/o1}

\end{document}